\documentclass[a4paper,12pt]{article}

\usepackage{listings}
\usepackage{longtable}

\usepackage{graphicx}

\usepackage[font=scriptsize]{caption}

\usepackage{float}

\usepackage{amsmath,amssymb,mathrsfs,eucal}

\usepackage[T1,T2A]{fontenc}

\font\cyrfont=wncyss10

\def\sza{\hbox{\cyrfont X}} 

\newtheorem{thm}{Theorem}

\newtheorem{conj}[thm]{Conjecture} 

\newtheorem{rem}[thm]{Remark}

\begin{document}

\title{ Elliptic curves with exceptionally large analytic order 
of the Tate-Shafarevich groups }

\author{Andrzej D\k{a}browski and Lucjan Szymaszkiewicz} 

\date{}

\maketitle{}

{\it Abstract}. We exhibit $88$ examples of rank zero elliptic curves 
over the rationals 
with $|\sza(E)| > 63408^2$, which was the largest previously known 
value for any explicit curve. Our record is an elliptic curve $E$ with 
$|\sza(E)| = 1029212^2 = 2^4\cdot 79^2 \cdot 3257^2$.  
We can use deep results by Kolyvagin, 
Kato, Skinner-Urban and Skinner to 
prove that, in some cases, these orders are the true orders of $\sza$.
For instance, $410536^2$ is the true order of $\sza(E)$ 
for $E= E_4(21,-233)$ from the table in section 2.3.



\section{Introduction}

Let $E$ be an elliptic curve defined over $\Bbb Q$ of conductor $N_E$, 
and let $L(E,s)$ denote its $L$-series. 
Let $\sza(E)$ be the Tate-Shafarevich group of $E$, 
$E(\Bbb Q)$ the group of rational points, and $R(E)$ the regulator,  
with respect to the N\'eron-Tate height pairing. 
Finally, let $\Omega_E$ be the least positive real period of the N\'eron 
differential of a global minimal Weierstrass equation for 
$E$, and define $C_{\infty}(E)=\Omega_E$ or $2\Omega_E$ 
according as $E(\Bbb R)$ is connected or not, and let $C_{\text{fin}}(E)$ 
denote the product of the Tamagawa factors of $E$ at the bad primes. 
The Euler product defining $L(E,s)$ converges for $\text{Re}\,s>3/2$. 
The modularity conjecture, proven by Wiles-Taylor-Diamond-Breuil-Conrad,   
implies that $L(E,s)$ has an analytic continuation to an entire function. 
The Birch and Swinnerton-Dyer conjecture relates the arithmetic data 
of $E$ to the behaviour of $L(E,s)$ at $s=1$. 

\begin{conj} (Birch and Swinnerton-Dyer) (i) $L$-function $L(E,s)$ has a 
zero of order $r=\text{rank}\;E(\Bbb Q)$ at $s=1$, 

(ii) $\sza(E)$ is finite, and 
$$ 
\lim_{s\to1}\frac{L(E,s)}{(s-1)^r} =
\frac{C_{\infty}\, (E)C_{\text{fin}}(E)\, R(E)\, |\sza(E)|}{|E(\Bbb Q)_{\text{tors}}|^2}. 
$$
\end{conj} 
If $\sza(E)$ is finite, the work of Cassels and Tate shows that its order 
must be a square. 

The first general result in the direction of this conjecture was proven for 
elliptic curves $E$ with complex multiplication by  Coates and  Wiles in 1976 
\cite{CW}, who showed that if $L(E,1)\not =0$, then the group $E(\Bbb Q)$ is finite. 
Gross and  Zagier \cite{GZ} showed that if $L(E,s)$ has a first-order zero at $s=1$, 
then $E$ has a rational point of infinite order.  Rubin \cite{Rub} proves that 
if $E$ has complex multiplication and $L(E,1)\not =0$, then $\sza(E)$ is finite. 
Let $g_E$ be the rank of $E(\Bbb Q)$ and let $r_E$ the order of the zero 
of $L(E,s)$ at $s=1$. Then   Kolyvagin \cite{Kol} proved that, if $r_E\leq 1$, 
then $r_E=g_E$ and $\sza(E)$ is finite.  Very recently,  Bhargava,  Skinner 
and  Zhang \cite{BhSkZ} proved that at least $66.48 \%$ of all elliptic curves 
over $\Bbb Q$, when ordered by height, satisfy the weak form of the Birch 
and Swinnerton-Dyer conjecture, and have finite Tate-Shafarevich group.

When $E$ has complex multiplication by the ring of integers of an imaginary 
quadratic field $K$ and $L(E,1)$ is non-zero, the $p$-part of the Birch and 
Swinnerton-Dyer conjecture has been established by  Rubin \cite{Rub2} for 
all primes $p$ which do not divide the order of the group of roots of unity 
of $K$.  Coates et al. \cite{cltz} \cite{coa}, and  Gonzalez-Avil\'es 
\cite{G-A} showed that there is 
a large class of explicit quadratic twists of $X_0(49)$ whose 
complex $L$-series does not vanish at $s=1$, and for which the full 
Birch and Swinnerton-Dyer conjecture is valid (covering the case $p=2$ 
when $K=\Bbb Q(\sqrt{-7})$). 
The deep results by  Kato,   Skinner and  Urban (\cite{SkUr}, Theorem 2) 
and  Skinner (\cite{Sk}, Theorem C) 
allow, in specific cases (still assuming $L(E,1)$ is non-zero), to establish $p$-part of 
the Birch and Swinnerton-Dyer conjecture for elliptic curves  without complex multiplication 
for all odd primes $p$.

It has long been known that the order of $\sza(E)[p]$ can be arbitrarily 
large for elliptic curves $E$ defined over $\Bbb Q$ and $p=2, 3$ 
(for $p=3$, the result 
is due to  Cassels (\cite{Cas}), and for $p=2$ it is due to  McGuinness 
(\cite{McG}). It was later extended for $p=5$ by  Fisher \cite{F}, and 
for $p=7$ and $13$ by  Matsuno \cite{M2}, but no similar result is known for $p=11$ 
or $p>13$. Let us mention, that they all used the fact that there are infinitely 
many elliptic curves defined over $\mathbb Q$ with rational $p$-isogenies.   
We also stress that it has not yet been proven that there exist elliptic 
curves $E$ defined over $\Bbb Q$ for which $\sza(E)[p]$ is non-zero for 
arbitrarily large primes $p$. 

In our earlier papers, we have investigated 
(see \cite{dw}, \cite{djs}, \cite{DSz0}, \cite{DSz}, \cite{DSz2}) 
some numerical examples of $E$ defined over $\Bbb Q$ for which $L(E, 1)$ 
is non-zero and the order of $\sza(E)$ is large.   

We extend these numerical 
results here in this paper, with the largest proved example of $\sza(E)$ 
having order $410536^2 = 2^6\cdot 7^2 \cdot 7331^2$.  
We exhibit $88$ examples of rank zero elliptic curves 
over the rationals 
with $|\sza(E)| > 63408^2$, which was the largest previously known 
value for any explicit curve. 
For some of these examples we can use deep results by Kolyvagin, 
Kato, Skinner-Urban and Skinner, to 
prove that these orders are the true orders of $\sza$.

Our idea was to use the family $E_i(n,p)$ from \cite{dw} (see section 2.1 below), 
within the bounds $20\leq n\leq 24$ 
and $0 < |p| \leq 5000$ (the calculations in \cite{dw} were focused 
on the pairs $(n,p)$ within the bounds $3\leq n\leq 19$ and $0 < |p| \leq 1000$).
First step was to find good candidates, i.e. the curves $E_i(n,p)$ with rank zero 
and $\max_i |\sza(E_i(n,p))| > 50000^2$.  The next step was to calculate 
$|\sza(E_i(n,p))|$ exactly for all these good candidates.  
In these steps, we computed (or estimated) the analytic orders of $\sza(E)$, using 
the approximations to $L(E,1)$.
The computations were performed using the computer package  
PARI/GP \cite{PARI}.   
The total running time for the various computational parts was about 9 months. 

\bigskip 

This research was supported in part by PLGrid Infrastructure.
Our computations were carried out in 2016 on the Prometheus
supercomputer via PLGrid infrastructure.   We also used the HPC cluster HAL9000 
and desktop computers Core(TM) 2 Quad Q8300 4GB/8GB, all located at the 
Department of Mathematics and Physics of Szczecin University.

\section{Results}

The previously largest value for $|\sza(E)|$ was $63408^2$, found by 
D\k{a}browski and Wodzicki \cite{dw}. In (\cite{djs}, section 5) we 
propose a candidate with $|\sza(E)| > 100000^2$.  
Below we present the results of our recent search for 
elliptic curves with exceptionally large analytic order of the 
Tate-Shafarevich groups. We exibit $88$  examples of rank zero 
elliptic curves with $|\sza(E)| > 63408^2$. Our record is 
an elliptic curve $E=E_2(23,-348)$ with $|\sza(E)| = 1029212^2$. 
Also  note that the prime $19861$ divides the orders of $\sza(E_i(22,304))$ - 
the largest (at the moment) prime dividing the order of $\sza(E)$ of 
an elliptic curve over $\Bbb Q$.

\subsection{Preliminaries} 

In this section we compute the analytic order of $\sza(E)$, i.e., the quantity 
$$
|\sza(E)|={L(E,1)\cdot |E(\Bbb Q)_{\text{tors}}|^2\over
C_{\infty}(E)C_{\text{fin}}(E)},
$$
for certain special curves of rank zero. We use the following approximation of
$L(E,1)$ 
$$
S_m=2\sum_{n=1}^m{a_n\over n}e^{-\frac{2\pi n}{\sqrt{N}}},
$$
which, for
$$
 m\geq {\sqrt{N}\over 2\pi}\left(2\log 2+k\log 10-\log(1-e^{-2\pi/\sqrt{N}})\right),
$$
differs from $L(E,1)$ by less than $10^{-k}$.

\smallskip
Consider (as in \cite{dw}) the family
$$
E_1(n,p):\quad y^2=x(x+p)(x+p-4\cdot 3^{2n+1}),
$$
with 
$(n,p)\in\Bbb N\times(\Bbb Z\setminus\{0\})$. Any member of the family admits
three isogenous (over $\Bbb Q$) curves $E_i(n,p)$ ($i=2,3,4$):

\bigskip 

$
E_2(n,p): \quad y^2=x^3+4(2\cdot 3^{2n+1}-p)x^2+16\cdot 3^{4n+2}x, 
$

$
E_3(n,p): \quad y^2=x^3+2(4\cdot 3^{2n+1}+p)x^2+(4\cdot 3^{2n+1}-p)^2x, 
$

$ 
E_4(n,p): \quad y^2=x^3+2(p-8\cdot 3^{2n+1})x^2+p^2x. 
$ 

\bigskip 
In our calculations, 
we focused on the pairs of integers $(n,p)$ within the bounds $20\leq n\leq 24$ 
and $0 < |p| \leq 5000$. Recall that the calculations in \cite{dw} were focused 
on the pairs $(n,p)$ within the bounds $3\leq n\leq 19$ and $0 < |p| \leq 1000$.

\bigskip 

The conductors, $L$-series and ranks of isogenous curves coincide, 
what may differ is the orders of $E(\Bbb Q)_{\text{tors}}$ and $\sza(E)$,
the real period $\Omega_E$, and the Tamagawa number $C_{\text{fin}}(E)$. 
In our situation we are dealing with $2$-isogenies, thus the  analytic order 
of $\sza(E)$ can only change by a power of $2$.

\bigskip 
{\it Notation}. Let $N(n,p)$ denote the conductor of the curve $E_i(n,p)$. 
We put $|\sza_i|=|\sza(E_i)|$.

\subsection{Elliptic curves $E_i(n,p)$ with 
$50000^2 \leq \max(|\sza_i|) < 250000^2$}

\begin{center}
\footnotesize
\begin{longtable}{|r|r|r|r|r|r|r|r|} 
\hline 
\multicolumn{1}{|c|}{$(n,p)$} 
& 
\multicolumn{1}{|c|}{$N(n,p)$} 
& 
\multicolumn{1}{|c|}{$|\sza_1|$} 
& 
\multicolumn{1}{|c|}{$|\sza_2|$} 
& 
\multicolumn{1}{|c|}{$|\sza_3|$} 
& 
\multicolumn{1}{|c|}{$|\sza_4|$} 
\\ 
\hline 
\endhead 
\hline 
\multicolumn{6}{|r|}{{Continued on next page}} \\
\hline
\endfoot

\hline 
\hline
\endlastfoot 

$(20,-756)$ & 42551829106699251024 & $27993^2$ &  $55986^2$ & $27993^2$  & $27993^2$  \\ 

$(20,-2000)$ & 190293894141760627320 & $15081^2$ & $60324^2$  &  $15081^2$ &  $60324^2$ \\

$(20,192)$ & 109418989131512359065 & $3780^2$ & $60480^2$  &  $945^2$ &  $60480^2$ \\ 

$(22,-692)$ & 11978814802342833513168 & $15194^2$ &  $30388^2$ &  $7597^2$ & $60776^2$  \\

$(21,-128)$  & 1969541804367222465954 &  $34234^2$  &  $68468^2$  & $34234^2$  & $68468^2$  \\ 

$(20,-180)$ & 60788327295284644080 & $20970^2$ & $41940^2$  & $10485^2$  &  $83880^2$ \\ 
 
$(21, 3)$  &   31512668869875559452120  &   $10962^2$  &  $43848^2$  &   $5481^2$  &   $87696^2$ \\ 

$(20,-2448)$ & 1653442502431742344680 &  $22028^2$ &  $88112^2$ &  $22028^2$ &  $88112^2$ \\ 

$(20,2704)$ & 11379574869677285146824 & $48538^2$ & $97076^2$  &  $97076^2$ &  $48538^2$ \\

$(21,12)$ & 281363114909603209392 & $12768^2$ & $102144^2$ & $3192^2$ &  $102144^2$ \\ 

$(20,-608)$ & 16631686347989878669080 & $25787^2$  & $103148^2$  & $51574^2$  & $51574^2$  \\

$(21,192)$ & 984770902183611232737 &  $54648^2$ & $109296^2$  &  $27324^2$ &  $109296^2$ \\ 

$(20, 4788)$  &  25871512096873143639456  &  $27745^2$  &  $110980^2$  &  $27745^2$  &  $27745^2$ \\ 

$(20,2680)$  &  23938195173261478962720 &  $14474^2$  &  $115792^2$  &  $14474^2$  &  $57896^2$  \\  

$(20,-801)$  &  34625031227394133415352  &  $29338^2$  &  $58676^2$  &  $29338^2$  &  $117352^2$  \\ 

$(22, 1344)$   &  62040566837567507664447  & $60930^2$  & $121860^2$ & $30465^2$ & $60930^2$   \\ 

$(20,-1436)$ & 1832369310703810488288 & $32455^2$ & $129820^2$  & $32455^2$  & $129820^2$  \\ 

$(20,4768)$ & 10032879618827902147272 & $16254^2$ & $130032^2$  & $8127^2$  & $65016^2$  \\  

$(21,-24)$  &  31512668869875559452768  &  $34092^2$  &    $68184^2$  &    $17046^2$  &  $136368^2$  \\  

$(20,-1376)$  &  37640132261240251922904  &  $70010^2$  &  $140020^2$  &  $140020^2$  &  $70010^2$  \\ 

$(22,64)$ & 8862938119652501095881 & $72306^2$ & $144612^2$  & $36153^2$  & $144612^2$  \\ 

$(21,-1536)$ & 1969541804367222468066 & $75897^2$ &  $151794^2$ & $75897^2$  &  $151794^2$ \\ 
 
$(20,-6)$ & 14005630608833581979328 & $19248^2$ & $76992^2$  & $4812^2$  &  $153984^2$ \\ 

$(22,304)$  &  27493195799738370745848  &  $39722^2$  &  $158888^2$  &    $19861^2$  &  $79444^2$  \\  

$(21,1516)$ & 11663380372737145525968 & $25866^2$ &  $206928^2$ &  $12933^2$ &  $103464^2$ \\ 

$(21,480)$  &  39390836087344449300840  &  $54110^2$  &  $216440^2$  &  $27055^2$  &  $108220^2$  \\ 

$(23,1452)$ & 6451697601805864768272 & $55698^2$ & $222792^2$  & $27849^2$  &  $222792^2$ \\

\end{longtable}
\end{center}

\subsection{Elliptic curves $E_i(n,p)$ with 
$\max(|\sza_i|) \geq 250000^2$}

\begin{center}
\footnotesize
\begin{longtable}{|r|r|r|r|r|r|r|r|} 
\hline 
\multicolumn{1}{|c|}{$(n,p)$} 
& 
\multicolumn{1}{|c|}{$N(n,p)$} 
& 
\multicolumn{1}{|c|}{$|\sza_1|$} 
& 
\multicolumn{1}{|c|}{$|\sza_2|$} 
& 
\multicolumn{1}{|c|}{$|\sza_3|$} 
& 
\multicolumn{1}{|c|}{$|\sza_4|$} 
\\ 
\hline 
\endhead 
\hline 
\multicolumn{6}{|r|}{{Continued on next page}} \\
\hline
\endfoot

\hline 
\hline
\endlastfoot

$(21,4)$ & 15756334434937779726048 & $130614^2$ &  $261228^2$ & $65307^2$  & $261228^2$  \\ 

$(21,1248)$    &  102416173827095568122280  &  $70375^2$  &  $281500^2$  &  $70375^2$  &  $70375^2$  \\  

$(20,-201)$    &  234594312697962498467304 &  $141540^2$  &  $141540^2$  &  $283080^2$  &  $141540^2$   \\

$(23,960)$ &  398832215384362549313205 &  $96254^2$ & $385016^2$  & $48127^2$  & $192508^2$  \\ 

$(24,832)$ & 373306953599763346160205 & $75780^2$ & $303120^2$ & $37890^2$ & $151560^2$  \\   

$(20,1120)$ & 30637316956823460343320 & $20440^2$ &  $327040^2$ & $20440^2$  &  $81760^2$ \\ 

$(23,-84)$ & 17448909423065861532624 & $184991^2$ & $369982^2$  & $184991^2$  &  $184991^2$ \\

$(22,480)$  &  354517524786100043822760  & $99938^2$   & $399752^2$   & $49969^2$   & $199876^2$ \\ 

$(23,-8)$ & 7441767284139709375008  &  $102120^2$ & $204240^2$  & $51060^2$  &  $408480^2$ \\ 

$(21,-233)$   &  149845956054714394972728  &  $51317^2$  &  $205268^2$  &  $51317^2$ & $410536^2$ \\ 

$(23,-96)$  &  638131544614980078907464   &   $264696^2$  &  $529392^2$  &  $132348^2$  &  $529392^2$  \\ 

$(24,-96)$ & 302272836922885300534872  & $412146^2$  &  $824292^2$ &  $206073^2$ & $824292^2$  \\  

$(23,-348)$ & 37011629587668844576720608 & $514606^2$ & $1029212^2$ & $257303^2$ & $1029212^2$  \\ 

\end{longtable}
\end{center}

\subsection{Birch and Swinnerton-Dyer conjecture for elliptic curves with 
exceptionally large analytic order of Tate-Shafarevich groups}

In this subsection, we will use the deep results by Kato,  Skinner-Urban and Skinner, 
to prove the full version of the Birch-Swinnerton-Dyer conjecture for some  
elliptic curves $E_i(n,p)$ with exceptionally large analytic order 
of Tate-Shafarevich groups.  

Let $\overline{\rho}_{E,p}: 
\text{Gal}(\overline{\Bbb Q}/\Bbb Q) \to 
\text{GL}_2(\Bbb F_p)$ denote the Galois 
representation on the $p$-torsion of $E$. 
Assume $p\geq 3$.  
As mentioned in \cite{GJPST}, the work of  Kato \cite{Kato} together with the works of 
 Mazur and  Rubin \cite{MR},  and  Matsuno \cite{M}, implies the following result. 

\begin{thm} 
Let $E$ be an optimal elliptic curve over $\mathbb Q$ with conductor $N_E$. 
Assume that $p\not | 6N_E$, and $\overline{\rho}_{E,p}$ is surjective. 
If $L(E,1) \not= 0$, then $\sza(E)$ is finite and 
\begin{equation} 
\text{ord}_p(|\sza(E)|) \leq \text{ord}_p\left(
{ L(E,1)\over C_{\infty}(E)}\right). 
\end{equation}
\end{thm}

For elliptic curves without complex multiplication, 
the following theorem of  Skinner-Urban and Skinner implies that the $p$-part 
of the Birch and Swinnerton-Dyer conjecture holds for nice 
primes $p$ in the rank zero case.

\begin{thm} (\cite{SkUr}, Theorem 2; \cite{Sk}, Theorem C) Let $E$ be an 
elliptic curve over $\Bbb Q$ with conductor $N_E$. 
Suppose: (i) $E$ has good ordinary or multiplicative reduction at $p$; 
(ii) there exists a prime $q\not=p$ such that $q\mid\mid N_E$ 
and $\overline{\rho}_{E,p}$ is ramified at $q$; 
(iii) $\overline{\rho}_{E,p}$ is surjective.  
If moreover $L(E,1)\not=0$, then the $p$-part of the Birch 
and Swinnerton-Dyer conjecture holds true, that is,   
we have 
\begin{equation} 
\text{ord}_p(|\sza(E)|)=\text{ord}_p\left(
{|E(\Bbb Q)_{\text{tors}}|^2 L(E,1)\over C_{\infty}(E)C_{fin}(E)}\right). 
\end{equation} 
\end{thm} 

\begin{rem}  
(i) The condition (ii) can be removed in the good ordinary case by the results 
of X. Wan \cite{W}, so we will omit it in the calculations below.  
(ii) Of course, surjectivity of $\overline{\rho}_{E,p}$  implies its irreducibility. 
\end{rem}

For the following curves from the tables in 2.2 and 2.3, 
we can apply the above results by Kato,  Skinner-Urban and Skinner: 
$E_i(22,-692)$, $E_i(20,-608)$, $E_i(20,-1436)$, $E_i(20,4788)$, 
$E_i(22,304)$, $E_i(22,64)$, $E_i(21,1516)$, $E_i(21,4)$, $E_i(21,-233)$.

Let us give some details for the curves  $E_i=E_i(20,-1436)$. 
We can use the above results  
to show that $|\sza_1|=5^2 \times 6491^2$ is the true order of 
$\sza(E_1)$ (and, hence, all $|\sza_i|$ are the true orders 
of $\sza(E_i)$).  

\bigskip 

(i)  We have 
$N_{E_1} = 2^5 \times 3 \times 7 \times 31 \times 257 \times 359 \times 107323 \times 8883041$.

(ii) We have  
$j_{E_1}={2^6 285451^3 4660272567723053424015171049538317^3 
\over 3^{82} 7^8 31^2 257^2 359^2 107323^2 8883041^2}$. 
The representation  $\overline{\rho}_{E_1,p}$ is surjective 
for any prime $p\geq 19$ by Prop. 1.8 in \cite{Zyw}. 
On the other hand, Prop. 6.1 in \cite{Zyw}  
gives a criterion to determine whether $\overline{\rho}_{E,p}$ 
is surjective or not for any non-CM elliptic curve $E$ 
and any prime $p\leq 11$. For instance, the representation 
$\overline{\rho}_{E,5}$ is not surjective if and only if 
$j_E={5^3(t+1)(2t+1)^3(2t^2-3t+3)^3 \over (t^2+t-1)^5}$ or 
$j_E={5^2(t^2+10t+5)^3\over t^5}$ or $j_E=t^3(t^2+5t+40)$ 
for some $t\in\Bbb Q$. We have checked, using Pari/GP,  impossibility of these 
three cases to hold.   For the primes $p=13$ and $p=17$ (actually for any $p\geq 5$), 
we can give a less computational 
proof of surjectivity as follows (suggested by the referee).  First, it is known due to 
B. Mazur (see Theorem 3 in \cite{Ma} or Theorem 1.3 in \cite{Me}), 
that for any elliptic curve $E$ over 
$\Bbb Q$ with all its $2$-division points defined over $\Bbb Q$, the representation  
$\overline{\rho}_{E,p}$ is absolutely irreducible for any 
prime $p\geq 5$.  Using the irreducibility, and the fact that the denominator of 
$j_{E_1}$ is not a $p$th power, we can apply  
Lemmas 1 and 2 in Chapter IV, Section 3.2 of \cite{S}.

(iii) $E_1$ has good ordinary reduction at $5$ 
and $6491$: $(N_{E_1},5)=(N_{E_1},6491)$ $=1$, and $a_5(E_1)=2$, 
$a_{6491}(E_1)=108$.

(iv) $E_1$ has multiplicative reduction at any $p \in \{3, 7, 31, 257, 359, 107323,$ $  8883041\}$. 
Take $q=7$. Then $7 || N_{E_1}$, and 
$\overline{\rho}_{E_1,p}$ is ramified at $7$ for these $p\not= 7$, 
since these $p$'s do not divide $\text{ord}_7(\Delta_{E_1})$. 
Take $q=31$. Then $31 || N_{E_1}$, and 
$\overline{\rho}_{E_1,7}$ is ramified at $31$, 
since $7$ do not divide $\text{ord}_{31}(\Delta_{E_1})$.

(v)  Using Magma \cite{Magma} we see that $\sza(E_1)[2]$ is trivial, and it is easy to 
conclude that the order of $\sza(E_1)$ is odd.

\bigskip 

{\bf Remark.} 
The curves $E = E_i(23,-348)$ have additive reduction at $3$, and the above 
results by Kolyvagin, Kato, Skinner-Urban and Skinner are not enough to prove 
that the analytic orders of $\sza$ are the true ones in these cases.  
These methods only show that, say,  the true order of $\sza( E_2(23,-348))$ is 
$1029212^2 \cdot 3^{2k} = 2^4\cdot  3^{2k} \cdot 79^2 \cdot 3257^2$, 
for some non-negative integer $k$.  To have $k=0$, we additionally need 
to prove that $\sza(E)[3]$ is trivial.  One way is to use Magma \cite{Magma} to show 
that $\text{Sel}(E)[3]=0$.  We had run \texttt{ThreeSelmerGroup(E)} subroutine  
for about 2,5 months, but without success, and we stopped the calculations. 
 A second way is to use a result of Gross (see \cite{Gr}, Prop. 2.1, Prop. 2.3 or  
\cite{MN}, Theorem 0.5) and some Magma computations. 
Let $y_K \in E(K)$ be the basic Heegner point attached to $E$ and a suitable 
imaginary quadratic field $K$.  If $y_K \notin 3E(K)$, then  
we obtain $\sza(E/K)[3^{\infty}]=0$ as desired.  But we didn't manage to 
compute $y_K$, since all our attempts required recomputing $L(E,1)$.

\subsection{Values of the Goldfeld-Szpiro ratio}

Let 

$$ 
GS(E) := {|\sza(E)| \over \sqrt{N_E}} 
$$ 
denote the {\it Goldfeld-Szpiro ratio} of an elliptic curve $E$.   
Eleven examples of elliptic curves  
with $GS(E) \geq 1$ are given by de Weger \cite{We}, the largest value being $6.893...$ .  
Another forty-seven examples with  $GS(E) \geq 1$ are produced by Nitaj \cite{Ni}, 
his largest value of $GS(E)$ being $42.265 ...$ .  For all of these examples the 
conductor does not exceed $10^{10}$.  The article of 
D\k{a}browski and Wodzicki \cite{dw} produce two examples with 
$GS(E) \geq 1$  for curves with much larger conductors.

The largest values of $GS(E)$ that we observed for our curves are given below. 
The notation $E_{i,j}(n,p)$ means that the given values of $|\sza(E)|$ and 
$GS(E)$ are shared by the isogeneous curves $E_i(n,p)$ and $E_j(n,p)$.

\begin{center}
\footnotesize
\begin{longtable}{|r|r|r|r|r|} 
\hline 
\multicolumn{1}{|c|}{$E$} 
& 
\multicolumn{1}{|c|}{$|\sza(E)|$} 
& 
\multicolumn{1}{|c|}{$GS(E)$} 
\\ 
\hline 
\endhead 
\hline 
\multicolumn{3}{|r|}{{Continued on next page}} \\
\hline
\endfoot

\hline 
\hline
\endlastfoot

$E_4(23,-8)$ &  $408480^2$  & $1.9342096803...  $   \\ 

$E_{2,4}(24,96)$    &   $824292^2$  &  $1.2358410273...  $     \\  

$E_2(23,-84)$ &  $369982^2$  & $1.0362798350...  $   \\ 

$E_4(20,-180)$    &  $83880^2$   &  $0.9024159172...  $     \\  

$E_{2,4}(21,12)$ &  $102144^2$  & $0.6220025144...  $   \\ 

$E_{2,4}(23,1452)$    &   $222792^2$  &  $0.6179625870...  $     \\  

$E_2(20,1120)$ &  $327040^2$  & $0.6110494864...  $   \\ 

$E_{2,4}(21,4)$    &  $261228^2$  &  $0.5436405656...  $     \\  

$E_{2,4}(21,-1536)$    &  $151794^2$   &  $0.5191903468... $     \\

\end{longtable}
\end{center}

Note that our record elliptic curves $E = E_{2,4}(23,-348)$  have relatively small 
Goldfeld-Szpiro ratio:   $GS(E) = 0.1741167606...$ .

\bigskip 

{\bf Acknowledgement.} We thank an anonymous referee for the constructive 
criticism and comments which improved the final version.

Institute of Mathematics, University of Szczecin, Wielkopolska 15, 
70-451 Szczecin, Poland; E-mail addresses: andrzej.dabrowski@usz.edu.pl and 
dabrowskiandrzej7@gmail.com;  
lucjansz@gmail.com

\end{document}